\documentclass[12pt]{article}
\usepackage{amsfonts}
\usepackage{amsmath,latexsym,amssymb,amsfonts,amsbsy, amsthm}
\usepackage[usenames]{color}
\usepackage{xspace,colortbl}
\allowdisplaybreaks

\newcommand{\beq}{\begin{equation}}
\newcommand{\eeq}{\end{equation}}
\newcommand{\ben}{\begin{eqnarray}}
\newcommand{\een}{\end{eqnarray}}
\newcommand{\beno}{\begin{eqnarray*}}
\newcommand{\eeno}{\end{eqnarray*}}

\begin{document}

\title{\textbf{On a power-type coupled system of Monge-Amp\`{e}re equations}}
\author{Zhitao Zhang$^{a,}$\thanks{Corresponding author, ~supported by NSFC 11325107, 11271353, 11331010.}~~ Zexin Qi$^b$ \\
{\small $^a$ Academy of Mathematics and Systems Science,}\\
{\small the Chinese Academy of Sciences, Beijing 100190, P. R. China}\\
 {\small E-mail: zzt@math.ac.cn}\\
{\small $^b$ College of Mathematics and Information Science,}\\
{\small Henan Normal University, Xinxiang 453007, P. R. China}\\
{\small E-mail: qizedong@126.com} \ }
\date{}
\maketitle
\begin{abstract}
\noindent We study an elliptic system coupled by Monge-Amp\`{e}re equations:
 \begin{center}
  $\left\{
   \begin{array}{ll}
      det~D^{2}u_{1}={(-u_{2})}^\alpha, & \hbox{in  $\Omega,$} \\
      det~D^{2}u_{2}={(-u_{1})}^\beta, & \hbox{in $\Omega,$} \\
      u_{1}<0, u_{2}<0,& \hbox{in  $\Omega,$}\\
     u_{1}=u_{2}=0, & \hbox{on $ \partial \Omega,$}
   \end{array}
 \right.$
 \end{center}
here $\Omega$~is a smooth, bounded and strictly convex domain
in~$\mathbb{R}^{N}$,~$N\geq2,~\alpha >0,~\beta >0$. When $\Omega$ is
the unit ball in $\mathbb{R}^{N}$, we use index theory of fixed
points for completely continuous operators to get existence,
 uniqueness results and nonexistence of radial convex solutions under
some corresponding assumptions on $\alpha,\beta$. When $\alpha>0$,
$\beta>0$ and $\alpha\beta=N^2$
  we also study a corresponding eigenvalue problem in more general domains.
\end{abstract}


\renewcommand{\theequation}{\thesection.\arabic{equation}}
\setcounter{equation}{0}
\section{Introduction}
Consider the following system coupled by Monge-Amp\`{e}re equations:
\begin{equation}\label{1}
     \left\{
   \begin{array}{ll}
      det~D^{2}u_{1}={(-u_{2})}^\alpha, & \hbox{in  $\Omega,$} \\
      det~D^{2}u_{2}={(-u_{1})}^\beta, & \hbox{in  $\Omega,$} \\
      u_{1}<0, u_{2}<0,& \hbox{in  $\Omega,$}\\
     u_{1}=u_{2}=0, & \hbox{on $\partial \Omega.$}
   \end{array}
 \right.
\end{equation}
Here~$\Omega$~is a smooth, bounded and strictly convex domain
in~$\mathbb{R}^{N}$,~$N\geq2,~\alpha >0,~\beta >0$; $det~D^{2}u$
stands for the determinant of Hessian matrix $(\frac{\partial^2
u}{\partial x_i\partial x_j})$ of~$u$.

Monge-Amp\`{e}re equations are fully nonlinear second order PDEs,
and there are important applications in geometry and other
scientific fields.
Monge-Amp\`{e}re equations have been studied in the past years
\cite{6,20,26,35,42}. However, to our best knowledge, only a few
works have been devoted to coupled systems. We refer the reader to
\cite{27} where the author established a symmetry result for a
system, which arises in studying the relationship between two
noncompact convex surfaces in $\mathbb{R}^3$. It seems to be H. Wang
\cite{36}\cite{wh} who first considered systems for Monge-Amp\`{e}re
equations.
He investigated the following system of equations:
\begin{equation}\label{2}
  \left\{
  \begin{array}{ll}
 det~D^{2}u_{1}={f(-u_{2})}, & \hbox{in $  B,$ }\\
      det~D^{2}u_{2}={g(-u_{1})}, & \hbox{in $ B, $}\\
     u_{1}=u_{2}=0, & \hbox{on $ \partial B.$}
  \end{array}
\right.
\end{equation}
Here and in the following $B:=\{x\in\mathbb{R}^N: |x|<1\}$. By reducing it to a system coupled by ODEs and using the fixed point index, the author obtained the following results:\\
\\
\textbf{Theorem A}~(\cite{36}, Theorem~1.1)~~\textit{Suppose~$f,g:[0,\infty)\rightarrow[0,\infty)$~are continuous.}
\begin{description}
\item[(a)] \textit{If~$f_{0}=g_{0}=0$~and~$f_{\infty}=g_{\infty}=\infty$,~then~}(\ref{2})~\textit{has at least one nontrivial radial convex solution.}
\item[(b)] \textit{If~$f_{0}=g_{0}=\infty$~and~$f_{\infty}=g_{\infty}=0$,~then~}(\ref{2})~\textit{has at least one nontrivial radial convex solution.}\\
 \end{description}

The notations were
\begin{center}
$f_{0}:=\lim\limits_{x\rightarrow0^{+}}\dfrac{f(x)}{x^{N}},~~~~~~~~~~f_{\infty}:=\lim\limits_{x\rightarrow\infty}\dfrac{f(x)}{x^{N}}$.
\end{center}
The above theorem implies the solvability of~(\ref{2})~is related to
the asymptotic behavior of $f, g$ at zero and at infinity.
Obviously, it asserts the existence of a radial convex solution for
system (\ref{1}) if $\Omega=B$ and one of the following cases holds:
\begin{center}
 (1) $\alpha > N,~\beta > N $;\;\;\;\;\;\;\;\;\;\;\;
 (2) $\alpha < N,~\beta < N $.
\end{center}
\par
{\bf What we are curious about is, for the sublinear-superlinear
case, i.e. $\alpha<N,~\beta>N$,~does system (\ref{1}) admits a
radial convex solution when~$\Omega= B$ ?}
\par
We obtain that:\\
\\
\textbf{Theorem~1.1}~~\textit{Let~$\Omega = B$,~then~}(\ref{1})~\textit{has a radial convex solution if~$\alpha>0,~\beta>0$~and~$\alpha\beta\neq N^2$~.}\\
\\
\textbf{Theorem~1.2}~~\textit{Let~$\Omega=B,~\alpha>0,~\beta>0$~and~$\alpha\beta<N^2$,~then}~(\ref{1})~
\textit{has a unique radial convex solution.}\\
\\
\textbf{Theorem~1.3}~~\textit{Let~$\Omega=B,~\alpha>0,~\beta>0$~and~$\alpha\beta=N^2$,~then~}(\ref{1})\textit{~admits no radial convex solution.}
\\

We also give new existence results for the more general system
(\ref{2}) in Remark 2.1. Our main tool is the fixed point index in a
cone  used in \cite{36}. However, based on the idea of decoupling
method we will consider a composite operator. Besides, solutions in
our theorems are classical, see Remark 2.2.

As $\alpha\beta= N^{2}$, for the eigenvalue problem
\begin{equation}\label{3}
    \left\{
   \begin{array}{ll}
      det~D^{2}u_{1}=\lambda{(-u_{2})}^\alpha, & \hbox{in  $\Omega,$} \\
      det~D^{2}u_{2}=\mu{(-u_{1})}^\beta, & \hbox{in $\Omega,$} \\
      u_{1}<0, u_{2}<0,& \hbox{in  $\Omega,$}\\
     u_{1}=u_{2}=0, & \hbox{on $ \partial \Omega$,}
   \end{array}
 \right.\\
\end{equation}
with positive parameters
$\lambda$ and $\mu$,
we have:\\
\\
\textbf{Theorem~1.4}~~\textit{Suppose $\Omega\subset\mathbb{R}^{N}$ is a bounded, smooth and strictly convex domain. If $\alpha>0,~\beta>0$ and $\alpha\beta=N^2$, then system (\ref{3}) admits a convex solution if and only if $\lambda\mu^\frac{\alpha}{N}=C$, where $C$ is a positive constant depending on $N$, $\alpha$ and $\Omega$}.

We will use the decoupling technique  again  to prove the
assertion. The solution operator is chosen to be of abstract form
which will be specified in section 3. What's more, a generalized
Krein-Rutman theorem(\cite{23}) is used. As to regularity, by Theorem 1.1  and second paragraph of p1253 of \cite{34}), we
see any eigenvector (admissible weak solution) of (\ref{3})  belongs to $C^1(\overline{\Omega})\times C^1(\overline{\Omega})$.

Recall the eigenvalue problem of the Monge-Amp\`{e}re operator,
\begin{center}
  $\left\{
    \begin{array}{ll}
      det~D^2u =|\lambda u|^{N}, & \hbox{in $\Omega$,} \\
      u =0, & \hbox{on $\partial \Omega$.}
    \end{array}
  \right.$
\end{center}
In \cite{23,26,35}, the authors proved by different methods that the
above equation has a unique positive eigenvalue, called the
principal eigenvalue of the Monge-Amp\`{e}re operator. Now we
consider \begin{equation} \label{zzt}
  \left\{
    \begin{array}{ll}
      det~D^2u =|\lambda v|^{N}, & \hbox{in} \Omega, \\
      det~D^2v =|\lambda u|^{N}, & \hbox{in} \Omega, \\
      u =v=0, & \hbox{on} \partial \Omega.
    \end{array}
  \right.
\end{equation}
By Theorem 1.4, we immediately obtain the following result.\\
\\
\textbf{Corollary 1.5}~~\textit{The  system \eqref{zzt} admits
nontrivial solutions if and only if $|\lambda|=\lambda_1(\Omega)$,
where $\lambda_1(\Omega)$ is the principal eigenvalue of the
 Monge-Amp\`{e}re operator corresponding to $\Omega$.}\\

This paper is organized as follows. In section 2 we give the proofs
of Theorem 1.1-1.3. The eigenvalue problem (\ref{3}) is discussed in
section 3 and we prove Theorem 1.4 there.
\section{Results concerning radial solutions}
\numberwithin{equation}{section}
 \setcounter{equation}{0}
When $\Omega=B$, let us search radial convex classical ($C^2(\Omega)$) solutions of
(\ref{1}). One can convert it to the following system of ODEs (see
Appendix A.2 of \cite{17} or \cite{22}):
\begin{equation}\label{4}
    \left\{
       \begin{array}{ll}
       \left(\left( u_{1}' (t)\right)^{N}\right)' = Nt^{N-1}\left(-u_{2}(t)\right)^{\alpha}, & \hbox{ $0<t<1$,} \\
       \left(\left( u_{2}' (t)\right)^{N}\right)' = Nt^{N-1}\left(-u_{1}(t)\right)^{\beta}, & \hbox{ $0<t<1$,} \\
       u_{1}<0,u_{2}<0, & \hbox{$0\leq t<1$,}\\
         u_{1}' (0)=u_{2}' (0)=0, u_{1}(1)=u_{2}(1)=0.
       \end{array}
     \right.
\end{equation}\\
In fact, the conversion is reversible if we choose a suitable
working space. However, we would rather look for solutions of
(\ref{4}) in $C^1[0,1]\times C^1[0,1]$ first and discuss the
regularity later in Remark 2.2.  Solutions of  problem (\ref{4}) are
equivalent to fixed points of a certain operator, and we can tackle
more general systems. Equivalently, we seek positive concave
solutions for convenience by letting $v_1=-u_1,v_2=-u_2$, and we can
transform the above system to
\begin{equation}\label{5}
    \left\{
       \begin{array}{ll}
        \left(\left( -v_{1}' (t)\right)^{N}\right)' = Nt^{N-1}\left(v_{2}(t)\right)^{\alpha}, & \hbox{ $0<t<1$,} \\
        \left(\left( -v_{2}' (t)\right)^{N}\right)' = Nt^{N-1}\left(v_{1}(t)\right)^{\beta}, & \hbox{ $0<t<1$,} \\
        v_{1}>0,v_{2}>0, & \hbox{$0\leq t<1$,}\\
         v_{1}' (0)=v_{2}' (0)=0, v_{1}(1)=v_{2}(1)=0.
       \end{array}
     \right.
\end{equation}\\
Below we will keep most notations used in \cite{36}. Recall the following lemma about fixed point index in a cone.\\
\\
\textbf{Lemma~2.1}~(\cite{12})~~\textit{Let $E$ be a Banach space, $K$ a cone in $E$. For $r>0$, define $K_r=\{u\in K:~\|u\|<r\}$. Assume $T:~\overline{K_r}\rightarrow K$ is completely continuous, satisfying $Tx\neq x$, $\forall x\in \partial K_r = \{u\in K: \|u\|=r\}$}.
\begin{enumerate}
  \item \textit{If~$\|Tx\|\geq \|x\|$,~$\forall x\in \partial K_r$,~then $i (T, K_r, K) = 0$.}
  \item \textit{If~$\|Tx\|\leq \|x\|$,~$\forall x\in \partial K_r$,~then $i (T, K_r, K) = 1$.}
\end{enumerate}

Now take the Banach space  to be $C[0,1]:=X$ with supremum norm.~Let
$K\subset X$ be
\begin{center}
 $K := \{v\in X: v(t)\geq 0,t\in [0, 1],~\min\limits_{\frac{1}{4}\leq t\leq \frac{3}{4}}v(t)\geq\dfrac{1}{4}\|v\|\}$,
\end{center}
which is  a cone in $X$. Denote $K_r= \{u \in K: \|u\|<r\}$ as in
Lemma 2.1. We introduce two solution operators. For~$v\in
K$,~define~$T_i(i=1, 2): K\rightarrow X$ to be
\begin{center}
   $T_1 (v) (t) = \displaystyle\int_{t}^{1}\left(\displaystyle\int_{0}^{s}N\tau^{N-1}v^{\alpha}
(\tau)d\tau\right)^{\frac{1}{N}}ds,~\hbox{$t\in [0,1]$;}$
\end{center}
 \begin{center}
 $T_2 (v) (t) = \displaystyle\int_{t}^{1}\left(\displaystyle\int_{0}^{s}N\tau^{N-1}v^{\beta}
(\tau)d\tau\right)^{\frac{1}{N}}ds,~\hbox{$t\in [0,1]$.}$
 \end{center}
Note the image of each operator is a nonnegative concave
$C^{1}$-function on $[0, 1]$, so by Lemma 2.2 in \cite{36}, the
above two operators map $K$ into itself. Besides, both operators are
completely continuous by standard arguments.

Define a composite operator $T=T_1T_2$, which is also completely
continuous from $K$ to itself. Calculation shows that $(v_1, v_2)\in
C^1[0,1]\times C^1[0,1]$ solves (\ref{5}) if and only if $(v_1,
v_2)$ belongs to $ K\backslash\{0\}\times K\backslash\{0\}$ and
satisfies $v_1=T_1v_2, v_2=T_2v_1$. \par Thus if $v_1\in
K\backslash\{0\}$ is a fixed point of $T$, define $v_2=T_2v_1$, then
$v_2\in K\backslash\{0\}$ so that $(v_1, v_2)\in C^1[0,1]\times
C^1[0,1]$ solves (\ref{5}); conversely, if $(v_1, v_2)\in
C^1[0,1]\times C^1[0,1]$ solves (\ref{5}), then $v_1$ must be a
nonzero fixed point of $T$ in $K$. So our task is to search nonzero
fixed points of $T$. We are in a position to give the following
proof of Theorem 1.1.\\
 \par
 {\bf Proof of Theorem 1.1}. Let $\Gamma$ be the positive
number given by
\begin{equation}\label{0}
\Gamma = \displaystyle\int_{\frac{1}{4}}^{\frac{3}{4}}\left(\displaystyle\int_{\frac{1}{4}}^{s}N\tau^{N-1}
d\tau\right)^{\frac{1}{N}}ds.
\end{equation}
For each $v\in K$,
\begin{align}
 \|T_2 (v)\| &= \displaystyle\int_{0}^{1}\left(\displaystyle\int_{0}^{s}N\tau^{N-1}v^{\beta}(\tau)d\tau\right)^{\frac{1}{N}}ds\notag\\
             &\geq\displaystyle\int_{\frac{1}{4}}^{\frac{3}{4}}\left(\displaystyle\int_{\frac{1}{4}}^{s}N\tau^{N-1}v^{\beta}(\tau)d\tau\right)^{\frac{1}{N}}ds\notag\\
             &\geq\displaystyle\int_{\frac{1}{4}}^{\frac{3}{4}}\left(\displaystyle
             \int_{\frac{1}{4}}^{s}N\tau^{N-1}\left(\frac{1}{4}\|v\|\right)^{\beta}d\tau\right)^{\frac{1}{N}}ds\notag\\
             &=\Gamma\left(\frac{1}{4}\|v\|\right)^\frac{\beta}{N}.\notag
\end{align}
Similarly we obtain
\begin{center}
$\|T_1 (v)\| \geq \Gamma\left(\dfrac{1}{4}\|v\|\right)^\frac{\alpha}{N}.$
\end{center}
Hence
\begin{align}
    \|T(v)\| &=\|T_1T_2(v)\|\notag\\
             & \geq \Gamma\left(\dfrac{1}{4}\|T_2(v)\|\right)^\frac{\alpha}{N}\notag\\
             &\geq\Gamma\left(\dfrac{1}{4}\Gamma\left(\dfrac{1}{4}\|v\|\right)^\frac{\beta}{N}\right)^\frac{\alpha}{N},\notag
\end{align}
which yields
\begin{equation}\label{8}
     \|T(v)\| \geq \Gamma_{1}\|v\|^{\frac{\alpha\beta}{N^{2}}}.
\end{equation}
where $\Gamma_1$ is a positive number that depends on $\alpha,\beta$ and $N$.

On the other hand, for each $v\in K$,
\begin{align}
 \|T_2 (v)\| &= \displaystyle\int_{0}^{1}\left(\displaystyle\int_{0}^{s}N\tau^{N-1}v^{\beta}(\tau)d\tau\right)^{\frac{1}{N}}ds\notag\\
             &\leq\left(\displaystyle\int_{0}^{1}N\tau^{N-1}v^{\beta}(\tau)d\tau\right)^{\frac{1}{N}}\notag\\
             &\leq\left(\displaystyle\int_{0}^{1}N\tau^{N-1}\|v\|^{\beta}d\tau\right)^{\frac{1}{N}}\notag\\
             &=\|v\|^\frac{\beta}{N}.\notag
\end{align}
Similarly,
\begin{center}
   $\|T_1 (v)\|\leq\|v\|^{\frac{\alpha}{N}},$
\end{center}
thus
\begin{equation}\label{9}
    \|T(v)\|\leq\|T_2(v)\|^{\frac{\alpha}{N}}\leq\|v\|^{\frac{\alpha\beta}{N^{2}}}.
\end{equation}
We take into account the following two cases.
\begin{enumerate}
  \item $\alpha\beta > N^{2}.$\\
Choose $r_1$ such that $0<r_1<1$. For $v\in K$ satisfying
$\|v\|=r_1$, we have $\|Tv\| <\|v\|$ by (\ref{9}). On the other
hand, by the estimate (\ref{8}), we can take $r_2$ large such that
$r_2>r_1$, and for each $v\in K$ satisfying $\|v\|=r_2$ it holds
$\|Tv\|>\|v\|$. By Lemma 2.1,
\begin{center}
    $i(T , K_{r_1 }, K) = 1$,\;\;\;\;\;\;\;\;$i(T , K_{r_2 }, K) = 0$.
\end{center}
We obtain $i(T, K_{r_2}\backslash \overline{K_{r_1}}, K ) = -1$ due to the additivity of the fixed point index. Then by
the existence property of the fixed point index, $T$ has a fixed point say $v_1$ in $K_{r_2}\backslash \overline{K_{r_1}}$.
Denote $v_2 = T_2v_1$, then $(-v_1,-v_2)$ is the desired solution of (\ref{4}). Considering regularity(see Remark 2.2 below), we get a classical solution for system (\ref{1}) when $\Omega = B$.
  \item $\alpha\beta < N^{2}.$\\
By (\ref{8}), we can choose $r_3>0$ small enough such that for each
$v\in K$ satisfying $\|v\|=r_3$, it holds $\|Tv\|> \|v\|$. On the
other hand, the estimate (\ref{9}) ensures the existence of $r_4$
such that $r_4>r_3$ and for each $v\in K$ satisfying $\|v\|=r_4$, we
have $\|Tv\|<\|v\|$. By lemma 2.1, we get
\begin{center}
    $i (T , K_{r_4 }, K) = 1$,\;\;\;\;\;\;\;\;$ i (T , K_{r_3 }, K) = 0$.
\end{center}
The rest of the proof is similar to that in case 1 and we omit
it.\end{enumerate} \hfill{$\square$}
\\
\textbf{Remark 2.1}~~The right hand side of each equation in system
(\ref{4}) is of particular form, while we can handle more general
ones, i.e.
\begin{equation}\label{10}
    \left\{
       \begin{array}{ll}
       \left(\left( u_{1}' (t)\right)^{N}\right)' = Nt^{N-1}f(-u_2)(t), & \hbox{ $0<t<1$,} \\
       \left(\left( u_{2}' (t)\right)^{N}\right)' = Nt^{N-1}g(-u_1)(t), & \hbox{ $0<t<1$,} \\
       u_{1}<0,u_{2}<0, & \hbox{$0\leq t<1$,}\\
         u_{1}' (0)=u_{2}' (0)=0, u_{1}(1)=u_{2}(1)=0.
       \end{array}
     \right.
\end{equation}
Similar arguments go through and we can get the following
conclusion:\par
 If $f,g:[0,\infty)\rightarrow[0,\infty)$ are
continuous, both nondecreasing,
then (\ref{10}) admits a solution if one of the following cases is
satisfied.
\begin{enumerate}
  \item $\lim\limits_ {x \rightarrow 0^{+}} \dfrac{f^{\frac{1}{N}}\left(g^{\frac{1}{N}}(x)\right)}{x} = 0$ and $\lim\limits_ {x \rightarrow \infty} \dfrac{f^{\frac{1}{N}}\left(g^{\frac{1}{N}}(x)\right)}{x}=\infty;$
  \item $\lim\limits_{x \rightarrow \infty} \dfrac{f^{\frac{1}{N}}\left(g^{\frac{1}{N}}(x)\right)}{x} = 0$ and $\lim\limits_ {x \rightarrow 0^{+}} \dfrac{f^{\frac{1}{N}}\left(g^{\frac{1}{N}}(x)\right)}{x} = \infty.$
\end{enumerate}\hfill$\square$
\\
\\
\textbf{Remark~2.2}~~The solutions we obtained in Remark 2.1 are in $C^1[0,1]\times C^1[0,1]$. Suppose $(\overline{u}_1,\overline{u}_2)$ is a solution of system (\ref{10}), can we get classical solutions for system (\ref{2}) by letting $u_{1}(x)=\overline{u}_1(|x|), u_{2}(x)=\overline{u}_2(|x|)$ ? This is the case if $(\overline{u}_1,\overline{u}_2)$ has higher order regularity, say belongs to $(C^2[0,1)\cap C^1[0,1])\times (C^2[0,1)\cap C^1[0,1])$. To see this, we refer the reader to Lemma 3.1 of \cite{39}, which states that if $u(x)=\widetilde{u}(|x|)$ in B, then $u\in C^2(B)$ if and only if $\widetilde{u}\in C^2[0,1)$ and $\widetilde{u}'(0)=0$. So let us explore further the regularity of $(\overline{u}_1,\overline{u}_2)$. Since it is supposed a solution of system (\ref{10}), we have
    \begin{center}
    $\overline{u}_1(t)=-\displaystyle\int_{t}^{1}\left(\displaystyle\int_{0}^{s}N\tau^{N-1}f(-\overline{u}_2
(\tau))d\tau\right)^{\frac{1}{N}}ds,~\hbox{$t\in [0,1]$;}$
    \end{center}
\begin{center}
    $\overline{u}_1'(t)=\left(\displaystyle\int_{0}^{t}N\tau^{N-1}f(-\overline{u}_2
(\tau))d\tau\right)^{\frac{1}{N}},~\hbox{$t\in [0,1]$;}$
\end{center}
and
\begin{equation}\label{11}
    \overline{u}_1''(t)=\frac{1}{N}\left(\int_{0}^{t}N\tau^{N-1}f(-\overline{u}_2
(\tau))d\tau\right)^{\frac{1}{N}-1}(Nt^{N-1}f(-\overline{u}_2
(t))).
\end{equation}
Similarly, we can  obtain
\begin{equation}\label{12}
    \overline{u}_2''(t)=\frac{1}{N}\left(\int_{0}^{t}N\tau^{N-1}g(-\overline{u}_1
(\tau))d\tau\right)^{\frac{1}{N}-1}(Nt^{N-1}g(-\overline{u}_1
(t))).
\end{equation}
By (\ref{11}) and (\ref{12}), if $f(x)>0$ and $g(x)>0$ for arbitrary
$x>0$, then calculation shows $(\overline{u}_1,\overline{u}_2)$
belongs to $C^2[0,1]\times C^2[0,1]$. Thus we get a nontrivial
convex classical solution of system (\ref{2}), by letting
$u_{1}(x)=\overline{u}_1(|x|), u_{2}(x)=\overline{u}_2(|x|)$ on
$\overline{B}$.\hfill$\square$
\\

We turn to the proof of the uniqueness result. Fix $\alpha>0,\beta>0$ such that $\alpha\beta<N^2$ in system (\ref{4}). We only need to show $T$ has at most one fixed point in $K$. With this in mind, we will give a sketch of the proof, since the rest of it's idea is similar to that used in \cite{22} where uniqueness for one single equation was established.\\
\\
\textbf{Definition 2.1} (\cite{GL}, or \cite{22} Definition
3.1)~~\textit{Let~$P$ be a cone from a real Banach space $Y$. With
some $u_0\in P$ positive, $A: P\rightarrow P$ is called
$u_0$-sublinear if}
\begin{description}
  \item[(a)] \textit{for any $x>0$, there exists $\theta_1>0,\theta_2>0$ such that $\theta_1u_0\leq Ax \leq \theta_2u_0$};
  \item[(b)] \textit{for any $\theta_1u_0\leq x \leq \theta_2u_0$ and $t\in (0,1)$, there always exists some $\eta=\eta(x,t)>0$ such that $A(tx)\geq(1+\eta)tAx$}.
\end{description}
\textbf{Lemma~2.2} (\cite{GL}, or \cite{22} Lemma 3.3)~~\textit{An increasing and $u_{0}$-sublinear operator A can have at most one positive fixed point.}\\

Now we choose the Banach space to be $Y=X=C[0,1]$ as before, but we work in a new cone $P:=\{v\in Y: v(t)\geq0,t\in[0,1]\}$. Since $K\subset P$, we only need to show that $T$ has at most one
fixed point in $P$.\\
\par
\textbf{Proof of Theorem 1.2}. It is readily seen that $T_1,T_2$ are
increasing operators with respect to the partial order induced by
$P$. So is $T=T_1T_2$. By Lemma 2.2, we only need to verify
that $T$ is $u_0$-sublinear for some $u_0$ positive in $Y$. Since $\alpha\beta < N^2$, we can assume $\alpha<N$ without loss of generality(otherwise consider the operator $\overline{T}:=T_2T_1$). Under this assumption, take $u_0=1-t$, then $T_1$ satisfies (a) of Definition 2.1, which is a consequence of Lemma 3.4 in \cite{22}. From this we know $T=T_1T_2$ also satisfies (a) of Definition 2.1. The proof is complete if $T$ satisfies (b) of Definition 2.1. To this end, let $\theta_1u_0\leq x\leq\theta_2u_0$, $\xi\in (0,1)$, then direct calculation give $T_2(\xi x)=\xi^\frac{\beta}{N}T_2(x), T_1(\xi x)=\xi^\frac{\alpha}{N}T_1(x)$. Thus $T(\xi x)=T_1(\xi^\frac{\beta}{N}T_2(x))=\xi^\frac{\alpha\beta }{N^{2}}T_1T_2(x)\geq (1+\eta)\xi Tx$ for some $\eta >0$. The last inequality holds because $\xi \in (0, 1)$ and $\alpha\beta < N^2$.\hfill$\square$\\

Finally in this section, we prove the nonexistence result.\par
\textbf{Proof of Theorem 1.3}. As analyzed previously, we only need
to show that $T$ has no positive  fixed point in $K$. For each $v\in
K$, we have
\begin{equation}\label{zzc}
\begin{array}{lll}
 \|T_2 (v)\| &= \displaystyle\int_{0}^{1}\left(\displaystyle\int_{0}^{s}N\tau^{N-1}v^{\beta}(\tau)d\tau\right)^{\frac{1}{N}}ds\\ 
             &\leq\left(\displaystyle\int_{0}^{1}N\tau^{N-1}v^{\beta}(\tau)d\tau\right)^{\frac{1}{N}}\\
             &\leq\left(\displaystyle\int_{0}^{1}N\tau^{N-1}\|v\|^{\beta}d\tau\right)^{\frac{1}{N}}\\ 
             &=\|v\|^\frac{\beta}{N}.
             \end{array}
\end{equation}
Assume that $T$ has a positive fixed point $v_0$ in $K$, then $v_0$
must be a concave function satisfying $v_0(1)=0$ and
$v_0(t)>0,t\in[0,1)$. Thus if we take $v=v_0$ in the above estimate,
 we know the last inequality in \eqref{zzc} must be strict. Thus $\|T_2 v_0\|<\|v_0 \|^\frac{\beta}{N}$. Similarly we
 have
$\|T_1 v\|\leq \|v \|^\frac{\alpha}{N}$,  $\forall v\in K$.
Therefore, we obtain that  $\|Tv_0
\|<\|v_0\|^\frac{\alpha\beta}{N^{2}}=\|v_0\|$. This contradicts the
assumption $v_0 = Tv_0$.\hfill$\square$

\section{A corresponding eigenvalue problem}
\numberwithin{equation}{section}
 \setcounter{equation}{0}
Checking the proof of theorem 1.3, we see the argument would not go through if the radius of the ball is larger than 1. This observation leads us to consider the following system
\begin{equation}\label{13}
    \left\{
   \begin{array}{ll}
      det~D^{2}u_{1}={(-u_{2})}^\alpha,& \hbox{ in $B_{R},$} \\
      det~D^{2}u_{2}={(-u_{1})}^\beta, & \hbox{in $B_{R},$} \\
      u_{1}<0,u_{2}<0, & \hbox{in $B_{R},$} \\
     u_{1}=u_{2}=0, & \hbox{on  $\partial B_{R}.$}
   \end{array}
 \right.\\
\end{equation}
Here $B_{R}$ denotes the ball of radius $R$ centered at zero,
$\alpha,\beta>0$ are such that $\alpha\beta=N^2$. By scaling, the
solvability of (\ref{13}) is equivalent to that of the following
problem:
\begin{equation}\label{14}
    \left\{
   \begin{array}{ll}
      det~D^{2}u_{1}=\lambda{(-u_{2})}^\alpha,& \hbox{ in $B,$} \\
      det~D^{2}u_{2}=\mu{(-u_{1})}^\beta, & \hbox{in $B,$} \\
      u_{1}<0,u_{2}<0, & \hbox{in $B,$} \\
     u_{1}=u_{2}=0, & \hbox{on}  \partial B,
   \end{array}
 \right.\\
\end{equation} where $\lambda$ and $\mu$ are positive parameters.
By Theorem 1.3, (\ref{14}) admits no radial convex solution when
$\lambda=\mu=1$. Further calculations show  that  if \eqref{14} has a
radial solution, then $\lambda,\mu$ should be in a suitable range. Indeed, let $X$ be
$C[0,1]$ and the cone $K$ as in section 2. Now we consider the new
operators $\widetilde{T_1}$, $\widetilde{T_2}$ defined as:
  \begin{center}
   $\widetilde{T_1} (v) (t) = \displaystyle\int_{t}^{1}\left(\displaystyle\int_{0}^{s}N\tau^{N-1}\lambda v^{\alpha}
(\tau)d\tau\right)^{\frac{1}{N}}ds, ~t\in [0,1],~ v\in K;$
  \end{center}
\begin{center}
   $\widetilde{T_2} (v) (t) = \displaystyle\int_{t}^{1}\left(\displaystyle\int_{0}^{s}N\tau^{N-1}\mu v^{\beta}
(\tau)d\tau\right)^{\frac{1}{N}}ds, ~t\in [0,1], ~v\in K.$
  \end{center}
We also define $\widetilde{T}:=\widetilde{T_{1}}\widetilde{T_{2}}$,
and we will investigate the fixed points of $\widetilde{T}$.\par
 Notice that
\begin{center}
    $\|\widetilde{T_2}(v)\|\leq\mu^{\frac{1}{N}}\|v\|^\frac{\beta}{N}$,\;\;\;\;$\|\widetilde{T_1}(v)\|\leq\lambda^{\frac{1}{N}}\|v\|^\frac{\alpha}{N},$
\end{center}
which yield
\begin{center}
$\|\widetilde{T}(v)\|\leq\lambda^{\frac{1}{N}} \|\widetilde{T_2} (v)\|^\frac{\alpha}{N}
\leq\lambda^{\frac{1}{N}}\mu^{\frac{\alpha}{N^{2}}}\|v\|.$
\end{center}
So if $v\neq0$ is a fixed point of $\widetilde{T}$, we have
necessarily $\lambda\mu^{\frac{\alpha}{N}}\geq1$, which implies
$\lambda\mu^{\frac{\alpha}{N}}$ can't be too small. On the other
hand, with $\Gamma$ defined in (\ref{0}), we have for each $v\in K$,
\begin{align}
 \|\widetilde{T_2}(v)\|&=\mu^{\frac{1}{N}}\displaystyle\int_{0}^{1}\left(\displaystyle\int_{0}^{s}N\tau^{N-1}v^{\beta}(\tau)d\tau\right)^{\frac{1}{N}}ds\notag\\
             &\geq\mu^{\frac{1}{N}}\displaystyle\int_{\frac{1}{4}}^{\frac{3}{4}}\left(\displaystyle\int_{\frac{1}{4}}^{s}N\tau^{N-1}v^{\beta}(\tau)d\tau\right)^{\frac{1}{N}}ds\notag\\
             &\geq\mu^{\frac{1}{N}}\displaystyle\int_{\frac{1}{4}}^{\frac{3}{4}}\left(\displaystyle
             \int_{\frac{1}{4}}^{s}N\tau^{N-1}\left(\frac{1}{4}\|v\|\right)^{\beta}d\tau\right)^{\frac{1}{N}}ds\notag\\
             &=\mu^{\frac{1}{N}}\Gamma\left(\frac{1}{4}\|v\|\right)^\frac{\beta}{N}.\notag
\end{align}
Similarly,
\begin{center}
$\|\widetilde{T_1}(v)\|\geq \lambda^{\frac{1}{N}}\Gamma\left(\dfrac{1}{4}\|v\|\right)^\frac{\alpha}{N},$
\end{center}
hence
\begin{align}
 \|\widetilde{T}(v)\| &= \|\widetilde{T_1}\widetilde{T_2}(v)\|\notag\\
                      &\geq \lambda^\frac{1}{N}\Gamma\left(\dfrac{1 }{4}\|\widetilde{T_2}(v)\|\right)^\frac{\alpha}{N}\notag\\
                      &\geq \lambda^\frac{1}{N}\Gamma\left(\dfrac{1 }{4}\right)^ \frac{\alpha}{N}\left(\mu^\frac{1}{N}\Gamma \left(\dfrac{1}{4}\|v\|\right)^\frac{\beta}{N}\right)^\frac{\alpha}{N}\notag\\
                      &= \lambda^\frac{1}{N}\mu^\frac{\alpha}{N^{2}}\left(\dfrac{1}{4}\Gamma\right)^{1+\frac{\alpha}{N}}\|v\|.\notag
\end{align}
So if $v\neq0$ is a fixed point of $\widetilde{T}$, we have
necessarily
\begin{center}
 $\lambda^\frac{1}{N}\mu^\frac{\alpha}{N^{2}}\left(\dfrac{1}{4}\Gamma\right)^{1+\frac{\alpha}{N}}\leq1,$
\end{center}
which implies $\lambda \mu^{\frac{\alpha}{N}}$ can't be too large.

 Is equation (\ref{14}) solvable for suitable $\lambda$ and $\mu$? The answer is positive and the domain need not even be symmetric, as asserted by Theorem 1.4. Our tool for (\ref{3}) is a
 generalized Krein-Rutman theorem developed in \cite{23}, where the author discussed eigenvalue problems for a broader class of fully nonlinear elliptic operators, including the Monge-Amp\`{e}re
 operator. \par
 Recall some concepts first(see \cite{23} for details). Let $E$ be a real Banach space with a cone $M\subset E$. The partial order induced by $M$ is written: $u\preceq v\Longleftrightarrow v-u \in M$.
 Let $A: E \rightarrow E$. $A$ is said to be homogeneous if it is positively homogeneous with degree 1. $A$ is monotone if it satisfies $x\preceq y\Rightarrow A(x)\preceq A(y)$.
 $A$ is called positive if $A(M)\subseteq M$. Finally, a positive operator $A:E \rightarrow E $ is called strong(relative to $M$), if for all $u,v\in Im(A)\cap M \backslash\{0\}$,
  there exist positive constants $\rho$ and $\tau$(which may depend on $u,v$), such that $u-\rho v\in M$ and $v-\tau u\in M$. The main content of the generalized Krein-Rutman theorem given
  in \cite{23} is as follows.
  \\
\\
\textbf{Lemma 3.1} (\cite{23}, Theorem 2.7)~~\textit{Let $E$ contain a cone $M$. Let $A: E \rightarrow E$ be a completely continuous operator with} $A|_M : M \rightarrow M$ \textit{homogeneous, monotone, and strong. Furthermore, assume that there exist nonzero elements $w, A(w)\in Im(A)\cap M$. Then there exists a constant $\lambda_0>0$ with the following properties}:
\begin{enumerate}
  \item \textit{There exists $u\in M\backslash\{0\}$, with $u=\lambda_0A(u)$};
  \item \textit{If $v\in M\backslash\{0\}$ and $\lambda>0$ such that $v=\lambda A(v)$, then $\lambda=\lambda_0$}.
\end{enumerate}

We also need the following lemmas to prove Theorem 1.4.\\
\\
By Theorem 1.1  and second paragraph of p1253 of \cite{34}, we have
\\
\\
\textbf{Lemma~3.2}~ (A special case of Trudinger \cite{34}, Theorem
1.1)~\textit{Let $\Omega$ be a strictly convex bounded domain in
$\mathbb{R}^{N}$, $\psi \in C(\overline\Omega)$ with $\psi\geq 0$,
$\phi \in C(\overline{\Omega})$. Then there exists a unique
admissible weak solution $u \in C^{1}(\overline\Omega)$ of the
equation}
 \begin{equation}\label{15}
 \left\{
   \begin{array}{ll}
      det~D^2u = \psi, & \hbox{ in $\Omega$,} \\
     u = \phi, & \hbox{ on $\partial \Omega$.}
    \end{array}
 \right.
  \end{equation}
  \par
   The definition of admissible weak solution coincides with the Aleksandrov sense weak solution (please see  page 1252-1253 in Trudinger \cite{34}), so
    Lemma 3.2 is valid for the Aleksandrov sense weak solution in the following Remark 3.1,
   we use the Aleksandrov sense weak solution here and in the following part of this paper.
\\
\\
\textbf{Remark 3.1}~The admissible weak solution in Lemma 3.2 can be
viewed as in Aleksandrov sense.
Recall the notion of Aleksandrov solution(see \cite{20}, Definition
1.1.1, Theorem 1.1.13 and Definition 1.2.1). Let
$\Omega\subset\mathbb{R}^{N}$ be an open subset and $u:
\Omega\rightarrow \mathbb{R}$. The normal mapping of $u$, or
subdifferential of $u$, is the set-valued function $\partial
u:\Omega\rightarrow 2^{\mathbb{R}^{N}}$ defined by
\begin{center}
    $\partial u(x_{0})=\{p\in\mathbb{R}^{N}:u(x)\geq u(x_{0})+p\cdot (x-x_{0}),~\forall x\in \Omega\}$.
\end{center}
Given $e\subset\Omega$, define $\partial u(e)=\bigcup_{x\in e}\partial u(x)$.

Let $u$ be continuous, then the class
\begin{center}
    $\mathcal{S}=\{e\subset\Omega:\partial u(e)\;is\;Lebesgue\;measurable\} $
\end{center}
is a Borel~$\sigma$-algebra. The set function $Mu:\mathcal{S}\rightarrow\overline{\mathbb{R}},~Mu(e)=|\partial u(e)|$ is a measure, finite on compacts, that is called the Monge-Amp\`{e}re measure associated with the function $u$.

Let $\nu$ be a Borel measure defined in $\Omega$, an open and convex
subset of $\mathbb{R}^{N}$. The convex function $u\in C(\Omega)$ is
called {\bf a generalized solution or Aleksandrov solution} to the
Monge-Amp\`{e}re equation
\begin{center}
    det $D^{2}u=\nu$
\end{center}
if the Monge-Amp\`{e}re measure $Mu$ associated with $u$ equals $\nu$.\hfill$\square$
\\
\\
\textbf{Lemma 3.3}~(Comparison Principle, \cite{20})~~\textit{Let
$\Omega$ be a bounded convex domain in $\mathbb{R}^{N}$.
 Denote $\mu[u]$ the} \textit{Monge-Amp\`{e}re} \textit{measure determined by $u$. Let $u,v \in C(\overline{\Omega})$ be two convex functions satisfying}
\begin{center}
  $\left\{
    \begin{array}{ll}
     \mu[u](e) \geq \mu[v](e), & \hbox{ $\forall$ Borel  $e\subset\Omega$;} \\
      u \leq v, & \hbox{on $\partial \Omega$.}
    \end{array}
  \right.$
\end{center}
\textit{then $u(x)\leq v(x)$ for any $x\in\Omega$.}\hfill$\square$
\\
\\
We are ready to give the proof of Theorem 1.4.\\
 \par {\bf Proof of Theorem 1.4}. Let $E$ be the Banach space $C(\overline{\Omega})$ with
supremum norm. Choose the negative cone
 $M :=\{u\in E:u(x)\leq0,\forall x\in\Omega\}$. Notice the partial order induced by $M$ reads: $u\preceq v \Longleftrightarrow v(x)\leq u(x),\forall
 x\in\Omega$. Define $A_1:E\rightarrow E, A_1(u)=v$, where $v$ is the unique admissible weak solution (Aleksandrov solution) of the equation
\begin{equation}\label{A1}
  \left\{
    \begin{array}{ll}
     det~D^2v = |u|^{\alpha}, & \hbox{in $\Omega$,} \\
     v = 0, & \hbox{on $\partial \Omega$.}
    \end{array}
  \right.
\end{equation}
By Lemma 3.2, $A_1$ is well defined. Similarly we define
$A_2:E\rightarrow E, A_2(u)=v$, where $v$ is the unique admissible
weak solution (Aleksandrov solution) of the equation
\begin{equation}\label{A2}
  \left\{
    \begin{array}{ll}
     det~D^2v = |u|^{\beta}, & \hbox{in $ \Omega$,} \\
     v = 0, & \hbox{on $ \partial \Omega$.}
    \end{array}
  \right.
\end{equation}
By Lemma 3.2
for the admissible weak solutions of \eqref{A1} and \eqref{A2}, we see $A_1u\in C^1(\overline{\Omega}), A_2u\in C^1(\overline{\Omega})$.
Finally we define a composite operator $A:=A_1A_2$.\par
 Let us verify $A$ satisfies the assumptions of Lemma 3.1.

Firstly, $A_{1}$, $A_{2}$ (thus $A$) are completely continuous by
Proposition 3.2 of \cite{23}. Since $A(E)\subseteq M$, $A$ is
positive. Let $t>0$, we have $A_2(tu)=t^{\frac{\beta}{N}}A_2(u),
A_1(tv)=t^{\frac{\alpha}{N}}A_1(u)$. As $\alpha\beta=N^2$, we deduce
$$A(tu)=A_1A_2(tu)=A_1(t^{\frac{\beta}{N}}A_2(u))=tA_1A_2(u)=tA(u),$$
which implies that $A$ is homogeneous. Besides, it is easy to get
that $A_1$, $A_2$ are monotone operators by Lemma 3.3, so is $A$.
\par To  see $A$ is strong, notice if $u\in Im(A)\cap M
\backslash\{0\}$, then there exists a $v\in E\backslash\{0\}$ such
that $u=A_1(A_2v)$. Now $A_2v$ is a nonzero convex function that is
strictly negative in $\Omega$,
by Lemma 3.2,
 we see $A_2v\in C^1(\overline{\Omega}), A_1(A_2v)\in C^1(\overline{\Omega})$.
  Then Lemma 3.4 of \cite{16} gives the exterior normal derivative satisfies
$u_\nu >0,\forall x\in \partial\Omega$, since $u$ is convex thus
subharmonic and $u(x)<0$ for $x\in\Omega$. Using these facts, one
can get by definition that $A$ is a strong operator.
\par
Finally, $\mathcal{N}(A)=\{0\}$ where $\mathcal{N}(A):=\{u\in
M:A(u)=0\}$. We see all assumptions in Lemma 3.1 are satisfied, then
there  exist $u_{\ast}\in M\backslash\{0\}$ and $\lambda_0>0$ such
that $u_{\ast}= \lambda_0A(u_{\ast})$. \par If we define
$v_{\ast}=A_2(u_{\ast})$, then $(u_{\ast},v_{\ast})$ must be a
solution of the following system
\begin{center}
   $\left\{
   \begin{array}{ll}
      det~D^{2}\left(\dfrac{u}{\lambda_0}\right) ={(-v)}^\alpha, & \hbox{ in $\Omega$}, \\
      det~D^{2}v = {(-u )}^\beta, & \hbox{ in $\Omega$}, \\
      u<0,v<0, & \hbox{in $\Omega$},\\
     u = v =0, & \hbox{ on $ \partial \Omega.$}
   \end{array}
 \right.$
   \end{center}
Furthermore, by the second conclusion of Lemma 3.1, if $u_1\in M\backslash\{0\}$ and $\lambda>0$ satisfy
 $u_1=\lambda A(u_1)$, then $\lambda=\lambda_0$. So the following system
\begin{equation}\label{16}
   \left\{
   \begin{array}{ll}
      det~D^{2}{u} = \widetilde{\lambda}{(-v)}^\alpha, & \hbox{ in $\Omega,$} \\
      det~D^{2}v = {(-u )}^\beta, & \hbox{ in $\Omega,$} \\
      u<0,v<0, & \hbox{in $\Omega,$}\\
     u = v =0, & \hbox{ on $ \partial \Omega$}
   \end{array}
 \right.
\end{equation}
admits a solution if and only if $\widetilde{\lambda}=\lambda_0^{N}$.

Now we show that (\ref{3}) has a convex solution if and only if
$\lambda \mu^{\frac{\alpha}{N}}=\lambda_0^{N}$, which implies the
first conclusion of Theorem 1.4. Indeed, if $(u, v)$ is a convex
solution of (\ref{3}), then from $det~D^2v = \mu (-u)^\beta$ we have
$det~D^2({\mu^{-\frac{1}{N}}}v)=(-u)^\beta$.\par
 Let
$\widetilde{v}={\mu^{-\frac{1}{N}}}v$, then $(-v )^\alpha =
\mu^\frac{\alpha }{N}(-\widetilde v )^\alpha$, and thus $det~D^2u
=\lambda(-v)^\alpha= \lambda\mu^\frac{\alpha }{N}(-\widetilde v
)^\alpha$. It is easily seen $(u,\widetilde v)$ is a convex solution
of (\ref{16}) if $\widetilde \lambda = \lambda\mu^\frac{\alpha
}{N}$. Since we have proved that (\ref{16}) admits a convex solution
only when $\widetilde{\lambda}=\lambda_0^{N}$, we get $\lambda
\mu^{\frac{\alpha}{N}}=\lambda_0^{N}$.\par
 On the other hand, assume
$\lambda \mu^{\frac{\alpha}{N}}=\lambda_0^{N}$, set
$\widetilde\lambda=\lambda\mu^\frac{\alpha}{N}$, then
$\widetilde\lambda=\lambda_0^{N}$ and (\ref{16}) admits a convex
solution, say $(u, v)$. Define $v^{\star}=\mu^\frac{1}{N}v$, then it
is easy to show that $(u ,v^{\star})$ is a convex solution of
(\ref{3}). This finishes the proof.
\hfill{$\square$}
\\
\\
{\bf Acknowledgement}. The authors thank the referees very much for their careful reading and useful suggestion.

\end{document}